\documentclass{amsart} 
\usepackage{amsmath}
\usepackage{amsmath} 
\usepackage{amsfonts, amssymb}
\usepackage{amsthm}
\usepackage[margin=3cm]{geometry}
\usepackage{hyperref}
\hypersetup{
    colorlinks=true,
    linkcolor=blue,
    citecolor=green,
    urlcolor=cyan
}
\newtheorem{theorem}{Theorem}[section] 
\newtheorem{corollary}[theorem]{Corollary} 
\newtheorem{definition}{Definition}[section] 
\newtheorem{remark}{Remark}[section] 
\subjclass[2020]{26A33; 44A20; 44A10; 42A38} 
\keywords{Fractional calculus;Integral transforms; Mellin transform; Fractional differential equations; Weighted fractional calculus Fractional calculus with respect to a function} 

\begin{document}

\title{The $\psi$-$\omega$-Mellin transform.}
\author{Gustavo Dorrego}
\address{Universidad Nacional del Nordeste, Corrientes, Argentina}
\email{gadorrego@exa.unne.edu.ar}

\author{Luciano Luque}
\address{Universidad Nacional del Nordeste, Corrientes, Argentina}
\email{lluque@exa.unne.edu.ar}

\author{Rubén Cerutti}
\address{Universidad Nacional del Nordeste, Corrientes, Argentina}
\email{rceruttiar@yahoo.com.ar}

\begin{abstract}
This manuscript introduces a generalization of the Mellin integral transform within the framework of weighted fractional calculus with respect to an increasing function. The proposed transform is much more suitable for working with fractional operators that involve a weight and are defined with respect to a function. This work also explores the connections between the Laplace and Fourier integral transforms in the same context. To achieve this, a new formulation of the weighted Fourier integral transform with respect to a function is presented, along with a new version of the bilateral Laplace transform. We study some of the properties of these new operators, obtain an inversion formula and a convolution theorem, and also present a practical application as an illustrative example.
\end{abstract}

\maketitle
\section{Introduction.}

       In recent decades, fractional calculus has evolved from a mathematical curiosity into a vibrant field of study, with integral transforms such as the Mellin, Laplace, and Fourier transforms serving as indispensable tools for solving fractional differential equations in a wide range of scientific and engineering disciplines (see, for example, \cite{Fernandez 3, Kiryakova} and references therein). Despite significant advancements, a consistently open area of research is the search for a unified theoretical framework to encompass the increasing diversity of fractional operators now present in the literature. A promising avenue towards this unification has been the introduction of generalized fractional operators, particularly "fractional calculus with respect to functions" ($\psi$-calculus) and "weighted fractional calculus" ($\omega$-calculus) \cite{Fernandez, Fernandez2, Agrawal 1}.

Recent literature has explored these two classes of operators individually. For instance, Aziz and Rehman \cite{Aziz} introduced a generalized Mellin transform with respect to a function $\psi(x)$, demonstrating its utility for fractional operators. Similarly, Fernandez et al. have rigorously investigated the underlying structure of weighted fractional calculus and fractional calculus with respect to functions, emphasising the power of "conjugation relations" (also known as "transmutation relations") to deduce new results from established classical theory \cite{Fernandez, Fernandez2}. This operational approach simplifies the study of complex operators and provides a powerful methodological tool for the field. However, a notable gap remains: a Mellin transform that fully unifies these two crucial generalizations---incorporating both the weight function $\omega(x)$ and the 'with respect to a function' parameter $\psi(x)$---has not been adequately developed. The generalized transform presented by Aziz et al. \cite{Aziz} is a special case of this broader class, as it lacks the explicit weight function $\omega(x)$.

Motivated by these prior contributions and the identified gap, this paper introduces a new generalized integral transform, which we term the $\omega$-$\psi$-Mellin transform. Our proposed transform is more general and intrinsically better suited for a wide class of weighted fractional operators defined with respect to a function. This work presents a systematic exposition of the basic properties of this new operator, including a new formulation for the weighted Fourier integral transform and a new version of the bilateral Laplace transform within the same context. We derive an inversion formula, establish a convolution theorem, and illustrate the practical applicability of our transform by solving a fractional differential equation.

The remainder of this paper is structured as follows. Section 2 presents preliminary definitions and results necessary for our development. In Section 3, we define the new $\omega$-$\psi$-Mellin transform and explore its fundamental properties. Section 4 presents the inversion formula and the convolution theorem. A practical application of our methodology to the solution of a fractional differential equation is provided in Section 5. Finally, we conclude with a summary of our findings and a discussion of potential future work in Section 6.

\section{Preliminaries.}
For the construction of our generalisation, it is fundamental to establish the theoretical framework and the key operators of weighted fractional calculus with respect to a function. This section presents definitions and theorems from previous works, principally by Aziz et al. \cite{Aziz} and Fernandez \cite{Fernandez}, which form the foundation for the development of this paper.

\begin{definition}\label{def:psi_mellin_transform}\cite{Aziz} 
Given $\psi\in C^{1}[a,b]$ an increasing function such that $\psi(0)=0$ and $\psi'(s)\neq 0$, $\forall s\in [a,b]$. The $\psi$-Mellin transform of a function $g$ is given by
\begin{equation}\label{eq:mellin_transform_def} 
        \mathcal{M}_{\psi}[g(x)](p)=G_{\psi}(p)=\int_{0}^{\infty}(\psi(x))^{p-1}g(x)\psi'(x)dx, \;\; p>0.
    \end{equation}
    While the inverse transform is given by
    \begin{equation}\label{eq:mellin_inverse_def} 
        \mathcal{M}^{-1}[G_{\psi}(p)]=\frac{1}{2\pi i}\int_{c-i\infty}^{c+i\infty}G_{\psi}(p)(\psi(x))^{-p}dp, \;\; p\in\mathbb{C}.
    \end{equation}
\end{definition}

\begin{definition}\label{def:psi_P_space}\cite{Aziz} 
    The space $${^{\psi}\mathcal{P}_a^b}:=\{(\psi(x))^{p-1}f(x)\in L_{\psi}^{1}(\mathbb{R}^+), a\leq \Re(s)\leq b\}$$ is a linear space for functions $h:\mathbb{R}^+\rightarrow\mathbb{C}$, where \\$L_{\psi}^{1}(\mathbb{R}^+)=\{f: f \;\text{is mesasurable on}\; \mathbb{R}^+, \int_0^{\infty}|f(s)|d(\psi)<\infty, \; d(\psi)=\psi' ds\}$.
\end{definition}

\begin{theorem}\label{thm:mellin_integral_convergence}\cite{Aziz} 
    Let $f\in{^\psi\mathcal{P}}^{b}_{a}$. Then, Mellin integral \eqref{eq:mellin_transform_def} converges absolutely and uniformly for $a\leq\Re(p)\leq b.$
\end{theorem}


\begin{definition}\cite{Fernandez}\label{def:w_weighted_RL_integral} 
	Let $\psi:[a,b]\rightarrow\mathbb{R}$ be a strictly increasing $C^1$ function, so that $\psi'>0$ everywhere, and let $w\in L^{\infty}(a,b)$ be a weigth function. The $w$-weighted Riemann-Liouville fractional integral with respect to $\psi$ of a given function $f\in L^1_{\psi}(a,b)$, to order $\alpha$ in $\mathbb{R}$ o $\mathbb{C}$, is defined by
	\begin{equation}\label{eq:RL_integral_formula} 
		^{RL}{_aI_{\psi,w}^{\alpha}}f(x)=\frac{1}{\Gamma(\alpha)w(x)}\int_{a}^{x}(\psi(x)-\psi(t))^{\alpha-1}w(t)f(t)\psi'(t)dt, \;\; x\in(a,b), 
	\end{equation} 
where $R(\alpha)>0$ if $\alpha\in\mathbb{C}$ or $\alpha>0$ if $\alpha\in\mathbb{R}$.
\end{definition}

Considering sufficient differentiability conditions for the functions $\psi$ and $\omega$, they can be defined:
\begin{itemize}
\item \textbf{Derivative in the Riemann-Liouville sense}(cf. \cite{Fernandez}) For $f\in AC_{\psi}^n[a,b]$ \begin{equation}\label{eq:RL_derivative_def} 
			{^{RL}_a D}_{\psi(x);w(x)}^{\alpha}f(x)=\left(\mathcal{D}_{\psi(x);w(x)}\right)^{n}{^{RL}_{a}I}_{\psi(x);w(x)}^{n-\alpha}f(x)
		\end{equation}
    \item \textbf{Derivative in the Caputo sense:}(cf. \cite{Fernandez}) For $f\in C_{\psi}^n[a,b]$
    \begin{equation}\label{eq:Caputo_derivative_def} 
			{^{C}_a D}_{\psi(x);w(x)}^{\alpha}f(x)={^{RL}{_{a}I}}_{\psi(x);w(x)}^{n-\alpha}\left(\mathcal{D}_{\psi(x);w(x)}\right)^{n}f(x)
		\end{equation}
    \item \textbf{Derivative in the Hilfer sense:} (cf.\cite{Fernandez2})
    \begin{equation}\label{eq:Hilfer_derivative_def} 
			{^{H}_a D}_{\psi(x);w(x)}^{\alpha,\beta}f(x)={^{RL}{_{a}I}}_{\psi(x);w(x)}^{\beta(n-\alpha)}\left(\mathcal{D}_{\psi(x);w(x)}\right)^{n}{^{RL}{_{a}I}}_{\psi(x);w(x)}^{(1-\beta)(n-\alpha)}f(x)
		\end{equation}
\end{itemize}
	 with $\Re(\alpha)\geq 0$, or $\alpha$ a real number such that $n-1\leq\Re(\alpha)<n$.

     \vspace{0,3cm}
     Here we consider the differential operator $\mathcal{D}$ defined as	\begin{equation}\label{eq:diferential_operator} 
			\mathcal{D}_{\psi(x);w(x)}f(x)=\frac{1}{w(x)\psi'(x)}\frac{d}{dx}\left(w(x)f(x)\right)
		\end{equation}
        or its equivalent: 
        \begin{equation}
           \mathcal{D}_{\psi(x);w(x)} =\frac{1}{\psi'(x)}\left(\frac{d}{dx}+\frac{\omega'(x)}{\omega(x)}\right)
        \end{equation}
Here it is important that the function $(\mathcal{D}_{\psi(x);w(x)})^kf(x)$ exists for $k=1,...,n$ and that for the case $k=n$ it is, respectively, continuous or $L^{1}_{\psi}$.
\subsection{Composition with operators $Q$ and $M$}
 
As can be seen in \cite{Fernandez}, using the following operators
	\begin{equation}\label{eq:operator_Q_psi} 
		Q_{\psi}f=f\circ\psi
	\end{equation}
	\begin{equation}\label{eq:operator_M_w} 
		(M_{w(x)}f)(x)=w(x)f(x)
	\end{equation}
    
    \vspace{0,5cm}
it is possible to write the Riemann-Liouville integral and fractional derivatives as follows:
	
	\begin{equation}\label{eq:RL_integral_composition} 
		{^{RL}_{a}I}_{\psi(x);w(x)}^{\alpha}=M^{-1}_{w(x)}\circ Q_{\psi}\circ \left({^{RL}_{\psi(a)}I}_{x}^{\alpha}\right)\circ Q^{-1}_{\psi}\circ M_{w(x)},
	\end{equation}
	
	\begin{equation}\label{eq:RL_derivative_composition} 
		{^{RL}_a D}_{\psi(x);w(x)}^{\alpha}=M^{-1}_{w(x)}\circ Q_{\psi}\circ \left({^{RL}_{\psi(a)}D}_{x}^{\alpha}\right)\circ Q^{-1}_{\psi}\circ M_{w(x)},
	\end{equation}

        \begin{equation}\label{eq:Caputo_derivative_composition} 
         {^{C}_a D}_{\psi(x);w(x)}^{\alpha}=M^{-1}_{w(x)}\circ Q_{\psi}\circ \left({^{C}_{\psi(a)}D}_{x}^{\alpha}\right)\circ Q^{-1}_{\psi}\circ M_{w(x)},   
        \end{equation}  

        \begin{equation}\label{eq:Hilfer_derivative_composition} 
         {^{H}_a D}_{\psi(x);w(x)}^{\alpha,\beta}=M^{-1}_{w(x)}\circ Q_{\psi}\circ \left({^{H}_{\psi(a)}D}_{x}^{\alpha,\beta}\right)\circ Q^{-1}_{\psi}\circ M_{w(x)}   
        \end{equation}
	i.e. it can be expressed by compositions in terms of the classical fractional operators.
\section{Results.}

This section introduces a generalized Mellin integral transform relevant to fractional calculus with weight and with respect to a function. The discussion covers its fundamental properties, inversion formula, and convolution theorem. We also investigate its connections to the generalized Laplace and Fourier transforms.

\subsection{Generalized Mellin transform.}

We begin by considering the generalised weighted Lebesgue space
$$L^{1}_{\omega,\psi}(\mathbf{R}^+)=\{f:\mathbf{R}^+\rightarrow\mathbf{C}\; \text{is measurable on}\;\mathbf{R}^+/ \|f\|_{L^{1}_{\omega,\psi}(\mathbf{R}^+)}<\infty\}$$ 
 	where
    $$\|f\|_{L^{1}_{\omega,\psi}(\mathbf{R}^+)
    }=\int_{\mathbf{R}^+}|\omega(x)f(x)|d(\psi),$$
    
and 

$$\psi'(x)dx=d(\psi).$$

Let $\psi:[0,\infty)\rightarrow\mathbf{R}$ be differentiable and such that $\psi'>0$, $\psi(0)=0$ and $\psi(x)\rightarrow\infty$ for $x\rightarrow\infty$; and let $\omega\in L^{\infty}(\mathbf{R}^+)$.

\begin{definition}\label{def:generalized_psi_omega_mellin} 
 Let $f\in L^1_{\omega,\psi}(\mathbf{R}^+)$,
    \begin{equation}\label{eq:generalized_mellin_integral} 
		\mathcal{M}_{\psi,\omega}[f(x)](p)=F_{\psi,\omega}(p)=\int_{0}^{\infty}(\psi(x))^{p-1}\omega(x)f(x)\psi'(x)dx, \;\; p\in\mathbb{C}.
	\end{equation}
\end{definition}	
Note that the convergence of the integral will usually occur (as in the classical case) in an open strip $\langle a,b\rangle$ such that $a<Re(p)<b$. The largest open strip $\langle a,b\rangle$ in which the integral converges is called the fundamental strip.

Indeed, if there exist $\alpha$ and $\beta$ such that the following conditions of asymptotic behaviour are satisfied:
\begin{itemize}
\item $\omega(x)f(x) = O((\psi(x))^\alpha)$ for $x \to 0^+$.
\item $\omega(x)f(x) = O((\psi(x))^\beta)$ for $x \to \infty$;
\end{itemize}
then we will have that the integral converges to the vertical strip of the complex plane defined by $-\alpha < \text{Re}(p) < -\beta$ with $\alpha>\beta$.

\begin{theorem}\label{thm:mellin_composition_form} 
\eqref{eq:generalized_mellin_integral} can also be expressed in terms of the operators \eqref{eq:operator_Q_psi} and \eqref{eq:operator_M_w} as follows:
	
\begin{equation}\label{eq:mellin_composition_formula} 
\mathcal{M}_{\psi,\omega}=\mathcal{M}\circ Q_{\psi}^{-1}\circ M_{w(x)}
\end{equation} where $\mathcal{M}$ is the classical Mellin transform.
\end{theorem}
Indeed:
\begin{eqnarray}
\left(Q_{\psi}^{-1} \circ M_{w(x)}\right)f(x)&=&Q_{\psi}^{-1}(M_{w(x)}f(x))\\
&=& Q_{\psi}^{-1}(w(x)f(x))\\
&=& w(\psi^{-1}(x))f(\psi^{-1}(x)).
\end{eqnarray}
Then, applying the classical Mellin transform and making the change $u = \psi^{-1}(x)$, we obtain
\begin{eqnarray}
\mathcal{M}[Q_{\psi}^{-1}(M_{w(x)}f(x))] &=& \int_0^{\infty} x^{p-1}w(\psi^{-1}(x))f(\psi^{-1}(x))dx\\
&=& \int_0^{\infty} (\psi(u))^{p-1}w(u)f(u)\psi'(u)du.
\end{eqnarray}

\begin{corollary}\label{cor:mellin_psi_case} 
If $\omega(x)=1$, then\begin{equation}\label{eq:mellin_psi_case_formula} 
\mathcal{M}_{\psi,1}[f(x)](p)=\int_{0}^{\infty}(\psi(x))^{p-1}f(x)\psi'(x)dx,
\end{equation} i.e. we obtain the Mellin transform introduced in \cite{Aziz} and given by \eqref{eq:mellin_transform_def}.
\end{corollary}

\begin{corollary}\label{cor:classic_mellin_case} 
If $\omega=1$ and $\psi(x)=x$, then the classical Mellin transform results
\begin{equation}\label{eq:classic_mellin_formula} 
\mathcal{M}_{x,1}[f(x)](p)=\int_{0}^{\infty}x^{p-1}f(x)dx,
\end{equation}
\end{corollary}
\subsection{Some properties.}

    Taking into account equation~\eqref{eq:mellin_composition_formula}, it is easy to prove the properties of the new transform in terms of the classical transform. 

Here are some important properties:

\begin{enumerate}
	\item (Shifting property) \begin{equation}\label{eq:prop_shift_mellin} 
		\mathcal{M}_{\psi,\omega}[(\psi(x))^{a}w(x)f(x)](p)=\mathcal{M}_{\psi,\omega}[f(x)](p+a), \;\; a\in\mathbb{C}.
	\end{equation}
    \item (Mellin transform of derivatives)

    Assuming that  $f$ and $\mathcal{D}^n_{\psi,\omega}f(x) \in L^1_{\omega,\psi}(\mathbf{R}^+)$ and such that 
$$\lim_{\substack{x\to\infty+\\ x\to 0+}}\omega(x)\left(\psi(x)\right)^{p-k-1}\mathcal{D}_{\psi(x),\omega(x)}^{n-k-1}f(x)=0$$ for all $k=0,1,...,n-1$, then   \begin{equation}\label{eq:prop_diff_mellin} 
		\mathcal{M}_{\psi,\omega}[\mathcal{D}^n_{\psi,\omega}f(x)](p)=\frac{\Gamma(1-p+n)}{\Gamma(1-p)}\mathcal{M}_{\psi,\omega}[f(x)](p-n)
	\end{equation}
    This property can be demonstrated using the corresponding result for the classical Mellin transform (cf.\cite{Kilbas} pp 21, f. (1.4.37)) the relation \eqref{eq:mellin_composition_formula} and the following equality
    \begin{equation}
        M_{\omega(x)}^{-1}\circ Q_{\psi(x)} \circ \frac{d}{dx} \circ Q_{\psi(x)}^{-1} \circ M_{\omega(x)}=\frac{1}{\psi'(x)}\left(\frac{d}{dx}+\frac{\omega'(x)}{\omega(x)}\right)=\mathcal{D}_{\psi(x),\omega(x)}
    \end{equation} proved in \cite{Fernandez} and whose version for the $n$-th derivative case is 
    \begin{equation}
         M_{\omega(x)}^{-1}\circ Q_{\psi(x)} \circ \left(\frac{d}{dx}\right)^n \circ Q_{\psi(x)}^{-1} \circ M_{\omega(x)}=\mathcal{D}^n_{\psi(x),\omega(x)}
    \end{equation}

    Indeed,
    \begin{eqnarray}\label{braquet}
\mathcal{M}_{\psi,\omega}\left[\mathcal{D}^n_{\psi,\omega}f(x)\right](p)&=&\left(\mathcal{M}\circ Q_\psi^{-1}\circ M_{\omega}\right)\circ \left[\left(M_{\omega}^{-1}\circ Q_{\psi}\circ\left(\frac{d}{dx}\right)^n\circ Q_{\psi}^{-1}\circ M_{\omega}\right)f(x)\right](p)\nonumber\\
        &=&\mathcal{M}\left[\left(\frac{d}{dx}\right)^n( Q_{\psi}^{-1}\circ M_{\omega}\circ f)(x)\right](p)\nonumber\\
        &=&\frac{\Gamma(1+n-p)}{\Gamma(1-p)}\mathcal{M}\left[( Q_{\psi}^{-1}\circ M_{\omega}\circ f)(x)\right](p-n) \nonumber\\
        &+&\sum_{k=0}^{n-1}\frac{\Gamma(1+k-p)}{\Gamma(1-p)}\left[x^{p-k-1}\left(\frac{d}{dx}\right)^{n-k-1}(Q_{\psi}^{-1}\circ M_\omega\circ f)(x)\right]_{0}^{\infty} \nonumber\\
        &=&\frac{\Gamma(1+n-p)}{\Gamma(1-p)}\mathcal{M}_{\psi,\omega}\left[f(x)\right](p-n)\nonumber\\
        &+&\sum_{k=0}^{n-1}\frac{\Gamma(1+k-p)}{\Gamma(1-p)}\left[x^{p-k-1}\left(\frac{d}{dx}\right)^{n-k-1}(Q_{\psi}^{-1}\circ M_\omega\circ f)(x)\right]_{0}^{\infty}\\
        \ \nonumber\\
        &=&\frac{\Gamma(1+n-p)}{\Gamma(1-p)}\mathcal{M}_{\psi,\omega}\left[f(x)\right](p-n) \nonumber
    \end{eqnarray}

    \begin{remark}    
      Note that the expression in brackets in \eqref{braquet} is equal to $0$ for all values of $k=0,1,...,n-1$ according to the hypothesis:
    \begin{eqnarray}
0&=&\lim_{\substack{x\to\infty+\\ x\to 0+}}[x^{p-k-1}\left(\frac{d}{dx}\right)^{n-k-1}(Q_{\psi}^{-1}\circ M_\omega\circ f)(x)]\\
       &=&\lim_{\substack{x\to\infty+\\ x\to 0+}}\omega(x)\left\{M_{\omega}^{-1}\circ Q_{\psi}\left[x^{p-k-1}\left(\frac{d}{dx}\right)^{n-k-1}\right]\circ Q_{\psi}^{-1}\circ M_\omega\circ f\right\}\\
       &=&\lim_{\substack{x\to\infty+\\ x\to 0+}}\omega(x)\left(\psi(x)\right)^{p-k-1}\left\{M_{\omega}^{-1}\circ Q_{\psi}\circ \left(\frac{d}{dx}\right)^{n-k-1}\circ Q_{\psi}^{-1}\circ M_\omega\circ f\right\}\\
       &=&\lim_{\substack{x\to\infty+\\ x\to 0+}}\omega(x)\left(\psi(x)\right)^{p-k-1}\mathcal{D}_{\psi(x),\omega(x)}^{n-k-1}f(x).
    \end{eqnarray}
    \end{remark}

	\item \textbf{(Mellin transform of Riemann-Liouville fractional integrals)}
    \begin{enumerate}
       \item If $\Re(\alpha)>0$, $\psi(x)^{p+\alpha-1}f(x)\in L^1_{\psi,\omega}(\mathbb{R^+})$ and $\Re(p)<1-\Re(\alpha),$ then \begin{equation}\label{eq:prop_RL_integral_mellin} 
		\mathcal{M}_{\psi,\omega}[{^{RL}{_{0}I}}_{\psi(x);w(x)}^{\alpha}f(x)](p)=\frac{\Gamma(1-p-\alpha)}{\Gamma(1-p)}\mathcal{M}_{\psi,\omega}[f(x)](p+\alpha), \;\;\Re(\alpha+p)<1.
	\end{equation}

        \item Let $-n 
< \Re(\alpha)\leq 1 - n, \; n = 1, 2, . . .$ $\psi(x)^{p+\alpha-1}f(x)\in L^1_{\psi,\omega}(\mathbb{R^+})$, $f\in C^{n}([0,b])$ $b$ being any positive number and $\Re(p)<1-\Re(\alpha),$ then \eqref{eq:prop_RL_integral_mellin} hold if 
\begin{equation}\label{eq: condicion de nulidad}
    \left(\psi(x)\right)^{p-k-1}({^{RL}I}_{\psi,\omega}^{\alpha-k-1}f(x))=0 \;\; \text{for}\; x=0, x=\infty, \;\; k=1,2,...,n.
\end{equation}
    \end{enumerate}
    To prove (a), we consider \eqref{eq:mellin_composition_formula}, \eqref{eq:RL_integral_composition} and the classical result of the Mellin transform of a Riemann-Liouville fractional integral. (See for example in \cite{Samko} Theorems 7.4 and 7.5)
 Indeed,
    \begin{eqnarray}
         \mathcal{M}_{\psi,\omega}[{^{RL}{_{0}I}}_{\psi(x);w(x)}^{\alpha}f(x)](p)&=& \mathcal{M}\circ Q_{\psi}^{-1}\circ M_{w(x)}\circ M^{-1}_{w(x)}\circ Q_{\psi}\circ \left({^{RL}_{\psi(0)}I}_{x}^{\alpha}\right)\circ Q^{-1}_{\psi}\circ M_{w(x)}\circ f(x)\nonumber\\
         &=&\mathcal{M}\circ  \left({^{RL}_{\psi(0)}I}_{x}^{\alpha}\right)\circ Q^{-1}_{\psi}\circ M_{w(x)}\circ f(x)\nonumber\\
         &=&\mathcal{M} \left[{^{RL}_{\psi(0)}I}_{x}^{\alpha}\left( Q^{-1}_{\psi}\circ M_{w(x)}\circ f(x)\right)\right](p)\nonumber\\
         &=&\frac{\Gamma(1-\alpha-p)}{\Gamma(1-p)}\mathcal{M} \left[ Q^{-1}_{\psi}\circ M_{w(x)}\circ f(x)\right](p+\alpha)\nonumber\\
         &=&\frac{\Gamma(1-\alpha-p)}{\Gamma(1-p)}\mathcal{M}_{\psi,\omega}[f(x)](p+\alpha)
    \end{eqnarray}

    To prove (b), we use the identity 
    \begin{equation}
{^{RL}{_{0}I}}_{\psi;w}^{\alpha}f(x)=\mathcal{D}_{\psi,\omega}^n {^{RL}{_{0}I}}_{\psi;w}^{\alpha+n}f(x).
    \end{equation}
Then, applying the generalised Mellin transform to both members and using part (a), it is obtained that
    
    \begin{eqnarray}
        \mathcal{M}_{\psi,\omega}[{^{RL}{_{0}I}}_{\psi;w}^{\alpha}f(x)](p)&=&\mathcal{M}_{\psi,\omega}[\mathcal{D}_{\psi,\omega}^n {^{RL}{_{0}I}}_{\psi;w}^{\alpha+n}f(x)](p)\nonumber\\
        &=&\frac{\Gamma(1+n-p)}{\Gamma(1-p)}\mathcal{M}_{\psi,\omega}[{^{RL}{_{0}I}}_{\psi;w}^{\alpha+n}f(x)](p-n)\nonumber\\
        &+&\sum_{k=0}^{n-1}\frac{\Gamma(1+k-p)}{\Gamma(1-p)}\left[\omega(x)(\psi(x))^{p-k-1}\mathcal{D}_{\psi,\omega}^{n-k-1}\left({^{RL}{_{0}I}}_{\psi;w}^{\alpha+n}f(x)\right)\right]^{\infty}_{0}\nonumber\\
        &=&\frac{\Gamma(1-p-\alpha)}{\Gamma(1-p)}\mathcal{M}_{\psi,\omega}[f(x)](p+\alpha)\nonumber\\
        &+&\sum_{k=0}^{n-1}\frac{\Gamma(1+k-p)}{\Gamma(1-p)}\left[\omega(x)(\psi(x))^{p-k-1}{^{RL}{_{0}I}}_{\psi;w}^{\alpha-k-1}f(x)\right]^{\infty}_{0}.
    \end{eqnarray}
    Finally, we obtain the desired result by considering \eqref{eq: condicion de nulidad}.
	\item \textbf{(Mellin transform of Riemann-Liouville fractional derivatives)}
    
Let $\Re(\alpha) > 0$, $n = [\Re(\alpha)] + 1$, $\psi(x)^{p-\alpha-1}f(x)\in L^1_{\psi,\omega}(\mathbb{R^+})$, $\Re(s) < 1 + \Re(\alpha)$ and the conditions
\begin{equation}
\lim_{x \to 0+} [\psi(x)^{s-k-1}(I^{n-\alpha}_{\psi(x),\omega(x)}f)(x)] = 0 \quad (k=0,1,\dots, n-1)
\end{equation}
and
\begin{equation}
\lim_{x \to \infty} [\psi(x)^{s-k-1}(I^{n-\alpha}_{\psi(x),\omega(x)}f)(x)] = 0 \quad (k=0,1,\dots, n-1)
\end{equation}
hold, then \begin{equation}\label{eq:prop_RL_deriv_mellin} 
		\mathcal{M}_{\psi,\omega}[{^{RL}D}_{\psi(x);w(x)}^{\alpha}f(x)](p)=\frac{\Gamma(1-p+\alpha)}{\Gamma(1-p)}\mathcal{M}_{\psi,\omega}[f(x)](p-\alpha)\quad (\Re(s-\alpha) < 1).
	\end{equation}
    To demonstrate this property, we use \eqref{eq:RL_derivative_composition}, \eqref{eq:mellin_composition_formula} and the classical result of the Mellin transform of the Riemann-Liouville derivative which can be seen for example in \cite{Kilbas}
.
 \item Let $\mu, \alpha\in\mathbf{C}$, $\Re(\alpha)>0$  \begin{equation}\label{eq:prop_mu_RL_deriv_mellin} 
 \mathcal{M}_{\psi,\omega}[(\psi(x))^{\mu}\omega(x){^{RL}D}_{\psi(x);w(x)}^{\alpha}f(x)](p)=\frac{\Gamma(1-\mu-p+\alpha)}{\Gamma(1-\mu-p)}\mathcal{M}_{\psi,\omega}[f(x)](\mu+p-\alpha)
	\end{equation}
    To demonstrate this property, property 1 and the previous property are used.
\vspace{1cm}

The following last two properties are still to be used \eqref{eq:mellin_composition_formula}, \eqref{eq:Caputo_derivative_composition} and \eqref{eq:Hilfer_derivative_composition}
       \item \begin{equation}\label{eq:prop_Caputo_deriv_mellin} 
		\mathcal{M}_{\psi,\omega}[{^{C}D}_{\psi(x);w(x)}^{\alpha}f(x)](p)=\frac{\Gamma(1-p+\alpha)}{\Gamma(1-p)}\mathcal{M}_{\psi,\omega}[f(x)](p-\alpha)
	\end{equation}
    \item \begin{equation}\label{eq:prop_Hilfer_deriv_mellin} 
		\mathcal{M}_{\psi,\omega}[{^{H}D}_{\psi(x);w(x)}^{\alpha,\beta}f(x)](p)=\frac{\Gamma(1-p+\alpha)}{\Gamma(1-p)}\mathcal{M}_{\psi,\omega}[f(x)](p-\alpha)
	\end{equation}
\end{enumerate}

\section{Inversion fórmula.}
Every integral transform has a corresponding inversion formula to recover the original function. We give an inversion formula for the transform defined in \eqref{eq:generalized_mellin_integral} that emerges intuitively from the operators in \eqref{eq:operator_Q_psi} and \eqref{eq:operator_M_w}.
	\begin{definition}\label{def:inverse_mellin_transform} 
	    \begin{equation}\label{eq:inverse_mellin_composition} 
 \mathcal{M}_{\psi,\omega}^{-1}=M_{w(x)}^{-1}\circ Q_{\psi}\circ \mathcal{M}^{-1}
	\end{equation}
    i.e. \begin{equation}\label{eq:inverse_mellin_formula} 
        \mathcal{M}_{\psi,\omega}^{-1}[F_{\psi,\omega}(p)](x) = \frac{1}{2\pi i  w(x)}\int_{\gamma-i\infty}^{\gamma+i\infty} F_{\psi,\omega}(p)(\psi(x))^{-p}dp.
    \end{equation}
    for all those functions for which its inverse transform exists in the classical sense.
	\end{definition}

Indeed, let us recall that the inverse Mellin transform is given by
\vspace{0.5cm}

\begin{equation}
    \mathcal{M}^{-1}[F(p)](x) = \frac{1}{2\pi i}\int_{\gamma-i\infty}^{\gamma+i\infty} F(p)x^{-p}dp \quad x>0, \gamma=\Re(p);
\end{equation}
\vspace{0.5cm}
Therefore, using first the operator $Q_\psi$ and then the operator $M^{-1}_{\omega}$ results in

$Q_{\psi}(\mathcal{M}^{-1}[F(p)](x)) = \mathcal{M}^{-1}[F(p)](\psi(x)) = \frac{1}{2\pi i}\int_{\gamma-i\infty}^{\gamma+i\infty} F(p)(\psi(x))^{-p}dp$

\vspace{0.5cm}

$M_{w(x)}^{-1} \circ Q_{\psi}(\mathcal{M}^{-1}[F(p)](x)) = \frac{1}{2\pi i  w(x)}\int_{\gamma-i\infty}^{\gamma+i\infty} F(p)(\psi(x))^{-p}dp$.
\vspace{0,5cm}

\vspace{0.5cm}
\subsection{Convolution.}
We now define a generalised convolution product and subsequently a generalised convolution theorem, observing that they generalise the classical results. 
\begin{definition}\label{def:generalized_convolution} 
	We define the convolution between two functions $f$ and $g$ to be
	\begin{equation}\label{eq:generalized_convolution_formula} 
	\left(f\ast_{\psi,w}g\right)(x)=\frac{1}{w(x)}\int_{0}^{\infty}w(s)f(s)w\left[\psi^{-1}\left(\frac{\psi(x)}{\psi(s)}\right)\right]g\left[\psi^{-1}\left(\frac{\psi(x)}{\psi(s)}\right)\right]\frac{\psi'(s)}{\psi(s)}ds.
	\end{equation}
\end{definition}
\vspace{0,2cm}
\begin{remark}\label{rem:convolution_aziz_case} 
Note that if $\omega(x)=1$ then the definition of convolution given by the equation \eqref{eq:generalized_convolution_formula} coincides with the definition given in \cite{Aziz}. And if, in addition, one takes $\psi(x)=x$ one obtains the classical Mellin definition of convolution.
\end{remark}

\begin{theorem}\label{thm:convolution_theorem} 
    If $f,g\in L^1_{\omega,\psi}(\mathbf{R}^+)$, the $f\ast_{\psi,w}g\in L^1_{\omega,\psi}(\mathbf{R}^+)$ and 
	\begin{equation}\label{eq:convolution_theorem_formula} 
	\mathcal{M}_{\psi,\omega}[\left(f\ast_{\psi,w}g\right)(x)](p)=\mathcal{M}_{\psi,\omega}[f(x)](p)\cdot \mathcal{M}_{\psi,\omega}[g(x)](p)
	\end{equation}
\end{theorem}

Indeed
\[
\mathcal{M}_{\psi, w}\left[\left(f\ast_{\psi,\omega}g\right)x\right](p) = \int_0^\infty \psi(x)^{p-1}\int_0^\infty w(s)f(s)w\left[\psi^{-1}\left(\frac{\psi(x)}{\psi(s)}\right)\right]g\left[\psi^{-1}\left(\frac{\psi(x)}{\psi(s)}\right)\right]\frac{\psi'(s)}{\psi(s)}\psi'(x)dsdx
\]

\[
= \int_0^\infty w(s)f(s)\left(\int_0^\infty (\psi(x))^{p-1}w\left[\psi^{-1}\left(\frac{\psi(x)}{\psi(s)}\right)\right]g\left[\psi^{-1}\left(\frac{\psi(x)}{\psi(s)}\right)\right]\psi'(x)dx\right)\frac{\psi'(s)}{\psi(s)}ds
\]
now, making the change of variable $t=\psi^{-1}\left(\frac{\psi(x)}{\psi(s)}\right)$ it has to $\psi(t)=\frac{\psi(x)}{\psi(s)}$ and therefore $\psi'(t)dt=\frac{\psi'(x)}{\psi'(s)}dx$ 
\[
= \int_0^\infty w(s)f(s)\left(\int_0^\infty (\psi(t)w(s))^{p-1}w(t)g(t)\psi'(t)\psi(s)dt\right)\frac{\psi'(s)}{\psi(s)}ds
\]

\[
= \int_0^\infty w(s)\psi(s)^{p-1}f(s)\psi'(s)ds \int_0^\infty (\psi(t))^{p-1}w(t)g(t)\psi'(t)dt
\]

\[=\mathcal{M}_{\psi,w}\left[f(s)\right](p) \cdot\mathcal{M}_{\psi,w}\left[g(t)\right](p).
\]

\begin{corollary}\label{cor:convolution_theorem_aziz} 
    If $\omega(x)=1$, then the result of the theorem coincides with the one presented in\cite{Aziz}.
\end{corollary}
\begin{corollary}\label{cor:convolution_theorem_classic} 

If we also consider $\psi(x)=x$ we obtain the well-known convolution theorem for the Mellin transform.
\end{corollary}

\subsection{Relationship with the Laplace and Fourier transforms.}

This section establishes relationships between the Mellin, Laplace, and Fourier transforms, demonstrating how they are interlinked via established change of variable formulae. We define the analogous weighted bilateral Laplace and Fourier transforms with respect to a given function."

\subsubsection{Relation to the Laplace transform with weight and with respect to a function.}
It should be noted that in the classic case, the following occurs
\begin{equation}\label{eq:laplace_mellin_classic_relation} 
M\{f(x)\}(s) = L_B\{f(e^{-t})\}(s)
\end{equation}
being
\begin{equation}\label{eq:LB_transform_def} 
L\{g(t)\}(s)+L\{g(-t)\}(s)=L_B\{g(t)\}(s)=\int_{-\infty}^{\infty}g(t)e^{-st}dt
\end{equation}
and by making the change $x=e^{-t}$ in the classical Mellin transform and using the relation \eqref{eq:mellin_composition_formula}, it can simply be proved that
\begin{equation}\label{eq:mellin_psi_omega_laplace_relation} 
\mathcal{M}_{\psi,\omega}[f(x)](p)={^B\mathcal{L}}_{\psi,\omega}[f(e^{-x})](p)
\end{equation}
 by defining the bilateral Laplace transform of a function with respect to another function and with a weight function $\omega$:
\begin{equation}\label{eq:bilateral_laplace_psi_omega_def} 
{^B\mathcal{L}}_{\psi,\omega}[f(x)](p)=\int_{-\infty}^{\infty}e^{-p\psi(x)}\omega(x)f(x)\psi'(x)dx
\end{equation}
\subsubsection{Relation to the Fourier transform with weight and with respect to a function.}
Similarly, if we define the Fourier transform of one function with respect to another and with a weight $\omega$ as follows
\begin{equation}\label{eq:fourier_psi_omega_def} 
    \mathcal{F}_{\psi,\omega}[f(x)](k)=\left(\mathcal{F}\circ Q_{\psi}^{-1}\circ M_{\omega}\right)[f(x)](k)=\frac{1}{\sqrt{2\pi}}\int_{-\infty}^{+\infty}e^{-ik\psi(x)}\omega(x)f(x)\psi'(x)dx
\end{equation}
and taking into account that in the classical case it is fulfilled that:
$$F\{g(t)\}(\omega) = M\{g(\ln x)\}(1 - i\omega),$$
we can show that the Mellin transform \eqref{eq:generalized_mellin_integral} relates to \eqref{eq:fourier_psi_omega_def} using the following formula:
\begin{equation}\label{eq:mellin_fourier_relation} 
    \mathcal{F}_{\psi,\omega}[f(x)](k)=\mathcal{M}_{\psi,\omega}[f(\ln x)](1-ip)
\end{equation}
See \cite{Fernandez tFourier} for Fourier transform of one function with respect to another.

\begin{remark}
    It is important to note that both the generalised bilateral Laplace transform \eqref{eq:bilateral_laplace_psi_omega_def} and the weighted Fourier transform with respect to a function \eqref{eq:fourier_psi_omega_def} are natural versions of these transforms that arise consistently within the context of the more general fractional calculus operators described in \cite{Fernandez}.
\end{remark}
\section{Application.}
To demonstrate the utility and power of the $\omega$-$\psi$-Mellin transform and the theorems presented herein, we apply our methodology to a fractional differential equation involving weighted Riemann-Liouville derivatives. This approach will be used to solve the following equation for $1<\alpha\leq 2$, subject to the boundary conditions 
$$\lim_{\substack{x\to\infty+\\ x\to 0+}}y(x)=0$$

\begin{equation}\label{eq:application_diff_eq} 
    (\psi(x))^{\alpha}\omega(x) {^{RL}D}_{\psi(x);w(x)}^{\alpha}y(x)=g(x).
\end{equation}

This class of equations is commonly found in the modelling of phenomena with memory effects, such as viscoelasticity, anomalous diffusion, or control systems.

Applying the generalised Mellin transform in both members
\begin{eqnarray}
\mathcal{M}_{\psi,\omega}\left[(\psi(x))^{\alpha}\omega(x) {^{RL}D}_{\psi(x);w(x)}^{\alpha}y(x)\right](p)&=&\mathcal{M}_{\psi,\omega}[g(x](p)\\
\frac{\Gamma(1-p)}{\Gamma(1-p-\alpha)}\mathcal{M}_{\psi,\omega}[y(x](p)&=&\mathcal{M}_{\psi,\omega}[g(x)](p)\\
\mathcal{M}_{\psi,\omega}[y(x](p)&=&\frac{\Gamma(1-p-\alpha)}{\Gamma(1-p)}\mathcal{M}_{\psi,\omega}[g(x)](p)
\end{eqnarray}

It turns out that
\begin{eqnarray}\label{eq:solution_transformed_y} 
    y(x)&=&\mathcal{M}^{-1}_{\psi,\omega}\left[\frac{\Gamma(1-p-\alpha)}{\Gamma(1-p)}G_{\psi,\omega}(p)\right](x)
\end{eqnarray}

The bracket of the second member can be interpreted as the product of Mellin transforms of the function $g$ and a function $h(x)=\mathcal{M}_{\psi,\omega}^{-1}\left[\frac{\Gamma(1-p-\alpha)}{\Gamma(1-p)}\right](x)$. Therefore, the solution of the equation is given by
\begin{equation}\label{eq:solution_convolution_form} 
    y(x)=(h\ast_{\psi,\omega} g)(x)
\end{equation}

But since we are looking for an explicit form for $h$, we use the relation \eqref{eq:inverse_mellin_composition} 
and the well-known result on the Mellin transform
\begin{equation}
    \mathcal{M}\left[\frac{(1-x)_{+}^{\alpha-1}}{\Gamma(\alpha)}\right](p)=\frac{\Gamma(1-\alpha-p)}{\Gamma(1-p)}, \;\; 0<\Re(\alpha)<1-\Re(p), \; (\textnormal{See~\cite{Samko},~p.,26}).
\end{equation}
\begin{eqnarray}\label{eq:explicit_h_function} 
h(x)=\mathcal{M}_{\psi,\omega}^{-1}\left[\frac{\Gamma(1-p-\alpha)}{\Gamma(1-p)}\right](x)&=&\left(M^{-1}_{\omega(x)}\circ Q_{\psi}\right)\left(\mathcal{M}^{-1}\left[\frac{\Gamma(1-p-\alpha)}{\Gamma(1-p)}\right](x)\right)\\
&=&M^{-1}_{\omega(x)}\left( Q_{\psi}\left(\frac{(1-x)_{+}^{\alpha-1}}{x^{\alpha}\Gamma(\alpha)}\right)\right)\\
&=&M^{-1}_{\omega(x)}\left(\frac{(1-\psi(x))_{+}^{\alpha-1}}{\psi(x)^\alpha\Gamma(\alpha)}\right)\\
&=&\frac{(1-\psi(x))_+^{\alpha-1}}{\omega(x)\psi(x)^\alpha\Gamma(\alpha)}.
\end{eqnarray}  

Recall that the function $t^{\alpha}_{+}$ is defined as 
\begin{equation}
    t^{\alpha}_+=
\begin{cases}
    t^\alpha       & \text{if } t>0;\\
    0      & \text{if } t\leq 0,
\end{cases}
\end{equation}

therefore it has to
\begin{equation}
    (1-\psi(x))_+^{\alpha-1}=
    \begin{cases}
    (1-\psi(x))^{\alpha-1}     & \text{if } 1>\psi(x);\\
    0      & \text{if } 1\leq \psi(x).
\end{cases}
\end{equation}

and returning now to the equation \eqref{eq:solution_convolution_form}, it turns out that the solution of the proposed problem is given by

\begin{eqnarray}\label{eq:final_solution_integration} 
    y(x)&=&\frac{1}{w(x)}\int_{0}^{\infty}w(s)h(s)w\left[\psi^{-1}\left(\frac{\psi(x)}{\psi(s)}\right)\right]g\left[\psi^{-1}\left(\frac{\psi(x)}{\psi(s)}\right)\right]\frac{\psi'(s)}{\psi(s)}ds\nonumber\\
    &=&\frac{1}{w(x)\Gamma(\alpha)}\int_{0}^{\psi^{-1}(1)}\frac{(1-\psi(s))_{+}^{\alpha-1}}{\psi(s)^\alpha}w\left[\psi^{-1}\left(\frac{\psi(x)}{\psi(s)}\right)\right]g\left[\psi^{-1}\left(\frac{\psi(x)}{\psi(s)}\right)\right]\frac{\psi'(s)}{\psi(s)}ds
\end{eqnarray}
\subsection{Particular cases of the solution.}

We now present three particular cases that can be derived from the solved general case.
\begin{itemize}
    \item Case 1: $\omega(x)=1$ and $\psi(x)=x$. 
    
    This case shows that the generalised transform reverts to the classical version.
    $$y(x) = \frac{1}{\Gamma(\alpha)} \int_0^1 (1-s)^{\alpha-1} s^{-\alpha-1} g\left(\frac{x}{s}\right) ds.$$

    \item Case 2: $\omega(x)=1$ and  $\psi(x)$ is an arbitrary increasing function.

       \begin{equation}
        y=\frac{1}{\Gamma(\alpha)}\int_{0}^{\psi^{-1}(1)}\frac{(1-\psi(s))_{+}^{\alpha-1}}{\psi(s)^\alpha}g\left[\psi^{-1}\left(\frac{\psi(x)}{\psi(s)}\right)\right]\frac{\psi'(s)}{\psi(s)}ds
    \end{equation}
    
    This solution complements the solution given in \cite{Aziz} in terms of another convolution product $\circ$ defined in that paper. Here we obtain the solution in terms of the convolution $\ast$.
    
    \item Case 3: $\omega(x)=1$, $\psi(x)=x^k$, $k>0$ and $g(x)=x^n$.

    Substituting these functions into the solution equation:
        $$y(x) = \frac{1}{\Gamma(\alpha)} \int_0^{1} (1-s^k)^{\alpha-1} (s^k)^{-\alpha} \left[\left(\frac{x^k}{s^k}\right)^{1/k}\right]^n \frac{k s^{k-1}}{s^k} ds,$$
        which can be reduced to
        $$y(x) = \frac{k x^n}{\Gamma(\alpha)} \int_0^1 (1-s^k)^{\alpha-1} s^{-\alpha k-n-1} ds.$$
   \item Case 4: $\omega(x) = e^x$, $\psi(x) = \ln(x+1)$, and $g(x) = x^2$, the solution is:

$$y(x) = \frac{1}{e^x\Gamma(\alpha)} \int_0^{e^1-1} \frac{(1-\ln(s+1))^{\alpha-1}_+}{(\ln(s+1))^\alpha}\left[\exp\left(e^{\frac{\ln(x+1)}{\ln(s+1)}}-1\right)\right]\\
\left[\left(e^{\frac{\ln(x+1)}{\ln(s+1)}}-1\right)^2\right] \frac{1}{s+1}\frac{ds}{\ln(s+1)}$$
    \item Case 5:  $\omega(x) = x^\alpha$ and $\psi(x) = \log(x)$, we obtain:
\[y(x) = \frac{1}{x^\alpha}\int_{0}^{\infty}s^\alpha h(s)\left(\psi^{-1}\left(\frac{\psi(x)}{\psi(s)}\right)\right)^\alpha g\left[\psi^{-1}\left(\frac{\psi(x)}{\psi(s)}\right)\right]\frac{1/s}{\log(s)}ds\]
where $\psi^{-1}(u)=e^u$, then $\psi^{-1}(\frac{\psi(x)}{\psi(s)})=e^{\frac{\log(x)}{\log(s)}}=x^{1/\log(s)}$. The solution, therefore, can be expressed as:
\[y(x) = \frac{1}{x^\alpha\Gamma(\alpha)}\int_{0}^{\infty}s^\alpha\frac{(1-\log(s))^{\alpha-1}_{+}}{s^\alpha}\left(x^{1/\log(s)}\right)^\alpha g\left(x^{1/\log(s)}\right)\frac{1/s}{\log(s)}ds\]
or equivalently:
\[y(x) = \frac{1}{x^\alpha\Gamma(\alpha)}\int_{0}^{\infty}(1-\log(s))^{\alpha-1}_{+}x^{\alpha/\log(s)}g\left(x^{1/\log(s)}\right)\frac{ds}{s\log(s)}.\]
\end{itemize}

  These particular cases demonstrate the versatility and generality of the proposed $\omega$-$\psi$-Mellin transform. By showing how the new transform encompasses existing transforms as special cases and provides solutions for complex, weighted fractional differential equations, this section validates its practical utility and theoretical coherence. It confirms that the new framework successfully unifies and extends prior work in the field, offering a robust tool for future research and application.

\section{Conclusions.}
This manuscript introduces and characterises a new generalisation of the Mellin integral transform, called the $\omega$-$\psi$-Mellin transform, within the framework of fractional calculus with a weight function $\omega$ and with respect to a function $\psi$. The motivation for this research was to fill a gap in the existing literature, as no transform had explicitly unified both generalisations in a single coherent framework.

The definition of the new operator has been formally established, it has been shown to encompass existing Mellin transforms as special cases, and its fundamental properties have been rigorously derived. The inversion formula and the convolution theorem, essential for application in solving equations, were developed with rigour. The utility and power of this new theoretical framework were practically illustrated by obtaining an analytical solution for a fractional differential equation of the Riemann-Liouville type.

Furthermore, the new formulations of the weighted and with respect to a function bilateral Laplace and Fourier transforms, presented in this work, represent an independent contribution that extends the scope of integral transform analysis into more general contexts.

This work lays a solid foundation for future research. The next directions could include the application of the Mellin transform to solve other types of fractional differential equations, as well as integral and integro-differential equations. Finally, exploring the properties of the new Laplace and Fourier transforms, as well as their applications in modelling physical phenomena with memory effects, represents a promising path for future investigations.

\section*{Statement.}
During the preparation of this work, the author(s) utilized Gemini 2.5 Flash with its Deep Research tool for language enhancement, interactive literature search and classification, improved search engine use, and coding assistance. Following the use of this tool/service, the author(s) reviewed and edited the content as necessary and take(s) full responsibility for the content of the published article.

\end{document}